\documentclass[a4paper]{article}

\usepackage[pages=all, color=black, position={current page.south}, placement=bottom, scale=1, opacity=1, vshift=5mm]{background}
\SetBgContents{
}      % copyright

\usepackage[linewidth=1pt]{mdframed}
\usepackage[margin=1in]{geometry}  
\usepackage{amsmath}
\usepackage{amsthm}
\usepackage{amssymb} 
\usepackage[utf8]{inputenc}
\usepackage{hyperref}
\hypersetup{
	unicode,
	pdfauthor={Author One, Author Two, Author Three},
	pdftitle={A simple article template},
	pdfsubject={A simple article template},
	pdfkeywords={article, template, simple},
	pdfproducer={LaTeX},
	pdfcreator={pdflatex}
}

\usepackage[T2A]{fontenc}   % For Cyrillic font encoding 

\usepackage[numbers,square]{natbib} 

\theoremstyle{plain}
\newtheorem{theorem}{Theorem}

\theoremstyle{definition}
\newtheorem{definition}[theorem]{Definition}
\newtheorem{example}[theorem]{Example}

\usepackage{graphicx, color}
\graphicspath{{fig/}}

\usepackage{algorithm, algpseudocode} % use algorithm and algorithmicx for typesetting algorithms
\usepackage{mathrsfs,ulem} % for \mathscr command

\title{Sarrus's Quilt: A Novel Way to Compute Determinants}
\author{%
    \begin{tabular}{c}
        \textsuperscript{1}Alexander D. Bonilla, \\
        \textsuperscript{2}Thomas Crawford \\
        \textsuperscript{3}Jorge Garcia \\
        \textsuperscript{4}Miles O'Brien \\
        \textsuperscript{5}Jasmine Torres \\[6pt]
        \textsuperscript{1,2,3,4,5}Mathematics Department, CSU Channel Islands, \\
        One University Drive, Camarillo, CA, USA \\[6pt]
        e-mail: 
        \textsuperscript{1}
        \href{mailto:alexander.bonilla871@csuci.edu}{\uline{alexander.bonilla871@csuci.edu}}, \\
        \textsuperscript{2} \href{mailto:thomas.crawford347@myci.csuci.edu}{\uline{thomas.crawford347@myci.csuci.edu}}\\   \textsuperscript{3}\href{mailto:jorge.garcia@csuci.edu}{\uline{jorge.garcia@csuci.edu}}\\ 
        \textsuperscript{4} \href{mailto:miles.obrien841@csuci.edu}{\uline{miles.obrien841@csuci.edu}}\\ 
        \textsuperscript{5} \href{mailto:jasmine.torres753@csuci.edu}{\uline{jasmine.torres753@csuci.edu}}
    \end{tabular}%
}
\date{} % Hapus tanggal otomatis

\usepackage{graphicx} % Required for inserting images
\usepackage[english]{babel}
\usepackage{amsthm}
\usepackage{xcolor}
\usepackage{tikz}
\usetikzlibrary{
    decorations.markings,
    decorations.pathmorphing,
    decorations.text,
}
\usetikzlibrary{calc,matrix}
\usepackage{amsmath}  

\usepackage{tkz-euclide}

\newcommand{\perm}[1]{$\left(\begin{smallmatrix} 1 & 2 & 3 & 4 & 5\\ #1 \end{smallmatrix}\right)$}

\newcommand{\Perm}[1]{\at{1 & 2 & 3 & 4 & 5\\ #1 }}

\usepackage{tikz-3dplot}

\newcommand{\minus}{\scalebox{0.50}[1.0]{$-$}}

\usepackage{tikz}
\usepackage{tikz-3dplot}   
\usetikzlibrary{calc,matrix}

\newcommand{\mat}[2][rrrrrrrrrrrrrrrrrrrrrrrrrrrrrrrrrrrrrrrrrr]{\left[
                     \begin{array}{#1}
                     #2\\
                     \end{array}
                    \right]}
                    
\newcommand{\at}[2][rrrrrrrrrrrrrrrrrrrrrrrrrrrrrrrrrrrrrrrrrr]{
                     \left(\begin{array}{#1}
                     #2\\
                     \end{array}\right)}
\newcommand{\pro}[5]{a_{#1}^1 a_{#2}^2 a_{#3}^3 a_{#4}^4 a_{#5}^5}

\newcommand{\xxmini}[2]{\begin{minipage}{#1\textwidth} \raggedleft #2\end{minipage}}

\newcommand{\keywords}[1]{\par\addvspace{17pt}\textbf{Keywords:}\ #1\par}

\begin{document}

\maketitle

\begin{abstract}
After analyzing the $4\times4$ determinant of a matrix, a shortcut was obtained to find such a determinant. Similarly to the Sarrus method for $2\times 2$ or $3\times 3$ determinants, the method consists of laying 19 columns of size 4 each and adding and subtracting some diagonal multiplications. There is a symmetry in the arrangement of these columns. A very symmetric pattern emerged for the $5\times5$ determinant which was also found. A cyclic pattern was observed for the diagonals of larger determinants when applying this method.
\end{abstract} 


\section{Introduction}
Pierre Frédéric Sarrus (1798–1861) was a French mathematician and professor who taught for over three decades at the University of Strasbourg. Although he contributed to several areas of mathematics, including the calculus of variations and mechanics, he is best remembered for his work in linear algebra and the invention of a straightforward method for computing the determinant of a $3
\times3$ matrix (Accademia delle Scienze, see~\cite{Accademia}).

Sarrus introduced his eponymous rule in 1833 as a compact method for evaluating the determinant of a $3\times3$ matrix. Prior to this, general determinant calculations relied on more elaborate cofactor expansion methods. The Rule of Sarrus provided a fast, mnemonic-based shortcut specific to the $3\times3$ case. It involves summing the products of diagonals running from upper left to lower right and subtracting the products of diagonals from lower left to upper right (Paul Cohn 1994, see~\cite{Cohn1994}). In Section~\ref{4x4}, we discuss a method for to compute the determinant of a $4\times4$ matrix that borrows from and extends the idea of the Sarrus rule. In Section~\ref{5x5} we compute the determinant of a $5\times5$ matrix that further expands on the Sarrus rule. In Section~\ref{LargerDeterm} we discuss the pattern that emerges in the signs of the diagonals of larger determinants and the cycle that appears in the pattern for the $n \times n$ case. Finally, in Section~\ref{ComparisonExisting} we examine other methods that have been published to compute the determinants of size $4\times4$, $5\times5$, and $n\times n$ matrices.

\section{The $4\times 4$ Determinant}\label{4x4}
We take a $4\times 4$ matrix and refer to its columns as columns 1, 2, 3, and 4 from left to right respectively. We will create a $4\times 19$ matrix using the columns of our original matrix. Similarly to how a $3\times 3$ determinant is solved with the Sarrus method, we write the columns 1-2-3-4, then repeat all but the last column, giving \textbf{1-2-3-4}-1-2-3. Following the Sarrus method, the diagonals give us some of the products produced within the determinant, but not all the ones we need. Additionally, diagonals starting from the first and third columns (from top left to bottom right and from bottom left to top right, giving a total of 4 diagonals) will be added to our "positive sum", whereas the diagonals produced from the second and fourth will be added to our "negative sum". By testing the products produced by shuffling the order of the columns, we deduced an order of the columns to give us all the necessary diagonals. The order of the columns will go \textbf{1-2-3-4}-1-2-\textbf{3-2-4-1}-3-2-\textbf{4-2-1-3}-4-2-1. Diagonals produced from columns 1, 3, 7, 9, 13, and 15 will be added together to create our "positive sum", while diagonals from columns 2, 4, 8, 10, 14, and 16 will be added together to create our "negative sum". For our final determinant, we simply have to subtract the negative sum from the positive sum.

We also noted that the pattern formed has a symmetry to it where if we fold it in half, 1 and 2 correspond to themselves, but 3 and 4 correspond to each other. \colorbox{green}{1-2}-3-4-\colorbox{green}{1-2}-3-\colorbox{green}{2}-4-\colorbox{green}{1}-3-\colorbox{green}{2}-4-\colorbox{green}{2-1}-3-4-\colorbox{green}{2-1}. In the following diagram, the blue diagonals will contribute to our "positive sum" regardless of the signs that each diagonal produces. Similarly, the orange diagonals contribute to our "negative sum".

 \ \,\hspace{0.25cm}\xxmini{0.5}{\begin{center}
\begin{tikzpicture}[
strip/.style = {
    draw=#1,%color
    line width=1em, opacity=0.2,
    shorten <=-2mm,shorten >=-2mm,
    line cap=round,
                            },
                    ]
\matrix (mtrx)  [matrix of math nodes,
                 column sep=0.5em,
                 nodes={text height=1ex,text width=2ex}]
{
|[blue]|+ & |[red]|- & |[blue]|+ & \color{blue}+\color{red}- & |[red]|- & |[blue]|+ & \color{blue}+\color{red}- & |[red]|- & |[blue]|+ & \color{blue}+\color{red}- & |[red]|- & |[blue]|+ & \color{blue}+\color{red}- & |[red]|- & |[blue]|+ & \color{blue}+\color{red}- & |[red]|- & |[blue]|+ & |[red]|-
                        \\[3.3 mm,between origins]
a_1 & b_1 & c_1 & d_1 & a_1 & b_1 & c_1 & b_1 & d_1 & a_1 & c_1 & b_1 & d_1 & b_1 & a_1 & c_1 & d_1 & b_1 & a_1\\
a_2 & b_2 & c_2 & d_2 & a_2 & b_2 & c_2 & b_2 & d_2 & a_2 & c_2 & b_2 & d_2 & b_2 & a_2 & c_2 & d_2 & b_2 & a_2\\
a_3 & b_3 & c_3 & d_3 & a_3 & b_3 & c_3 & b_3 & d_3 & a_3 & c_3 & b_3 & d_3 & b_3 & a_3 & c_3 & d_3 & b_3 & a_3\\
a_4 & b_4 & c_4 & d_4 & a_4 & b_4 & c_4 & b_4 & d_4 & a_4 & c_4 & b_4 & d_4 & b_4 & a_4 & c_4 & d_4 & b_4 & a_4\\[5 mm,between origins]
\colorbox{green}{1} & \colorbox{green}{2} & 3 & 4 & \colorbox{green}{1} & \colorbox{green}{2} & 3 & \colorbox{green}{2} & 4 & \colorbox{green}{1} & 3 & \colorbox{green}{2} & 4 & \colorbox{green}{2} & \colorbox{green}{1} & 3 & 4 & \colorbox{green}{2} & \colorbox{green}{1}\\
};
\draw[thick] (mtrx-2-1.north) -| (mtrx-5-1.south west)
                              -- (mtrx-5-1.south);
\draw[thick] (mtrx-2-4.north) -| (mtrx-5-4.south east)
                              -- (mtrx-5-4.south);
\path[draw,strip=orange]
    (mtrx-2-2.center) edge (mtrx-5-5.center)
    (mtrx-2-4.center) edge (mtrx-5-7.center)
    (mtrx-5-2.center) edge (mtrx-2-5.center)
    (mtrx-5-4.center) edge  (mtrx-2-7.center)
    (mtrx-2-8.center) edge (mtrx-5-11.center)
    (mtrx-2-10.center) edge (mtrx-5-13.center)
    (mtrx-5-8.center) edge (mtrx-2-11.center)
    (mtrx-5-10.center) edge  (mtrx-2-13.center)
    (mtrx-2-14.center) edge (mtrx-5-17.center)
    (mtrx-2-16.center) edge (mtrx-5-19.center)
    (mtrx-5-14.center) edge (mtrx-2-17.center)
    (mtrx-5-16.center)  --  (mtrx-2-19.center);
    
\path[draw,strip=cyan]
    (mtrx-2-1.center) edge (mtrx-5-4.center)
    (mtrx-2-3.center) edge (mtrx-5-6.center)
    (mtrx-5-1.center) edge (mtrx-2-4.center)
    (mtrx-5-3.center) edge (mtrx-2-6.center)
    (mtrx-2-7.center) edge (mtrx-5-10.center)
    (mtrx-2-9.center) edge (mtrx-5-12.center)
    (mtrx-5-7.center) edge (mtrx-2-10.center)
    (mtrx-5-9.center) edge (mtrx-2-12.center)
    (mtrx-2-13.center) edge (mtrx-5-16.center)
    (mtrx-2-15.center) edge (mtrx-5-18.center)
    (mtrx-5-13.center) edge (mtrx-2-16.center)
    (mtrx-5-15.center)  --  (mtrx-2-18.center);
\end{tikzpicture}
\end{center}}

\begin{example}  %Thomas Look Here!!!!!
  Let's compute the determinant of the following $4\times 4$ matrix, $\mat{2&3&4&-1\\ 1&-2&0&5\\ 5&2&2&-3\\ 8&1&1&1}$.
\newline Using the arrangement created above we now use our example matrix to make a new matrix based on that arrangement, shown below. 
\begin{center}
\begin{tikzpicture}[
strip/.style = {
    draw=#1,%color
    line width=1em, opacity=0.2,
    shorten <=-2mm,shorten >=-2mm,
    line cap=round,
                            },
                    ]
\matrix (mtrx)  [matrix of math nodes,
                 column sep=0.5em,
                 nodes={text height=1ex,text width=2ex}]
{
|[blue]|+ & |[red]|- & |[blue]|+ & \color{blue}+\color{red}- & |[red]|- & |[blue]|+ & \color{blue}+\color{red}- & |[red]|- & |[blue]|+ & \color{blue}+\color{red}- & |[red]|- & |[blue]|+ & \color{blue}+\color{red}- & |[red]|- & |[blue]|+ & \color{blue}+\color{red}- & |[red]|- & |[blue]|+ & |[red]|-
                        \\[3.3 mm,between origins]
2 &       3 & 4 & \minus1 & 2 &       3 & 4 &       3 & \minus1 & 2 & 4 &       3 & \minus1 &       3 & 2 & 4 & \minus1 &       3 & 2\\
1 & \minus2 & 0 &       5 & 1 & \minus2 & 0 & \minus2 &       5 & 1 & 0 & \minus2 &       5 & \minus2 & 1 & 0 &       5 & \minus2 & 1\\
5 &       2 & 2 & \minus3 & 5 &       2 & 2 &       2 & \minus3 & 5 & 2 &       2 & \minus3 &       2 & 5 & 2 & \minus3 &       2 & 5\\
8 &       1 & 1 &       1 & 8 &       1 & 1 &       1 &       1 & 8 & 1 &       1 &       1 &       1 & 8 & 1 &       1 &       1 & 8\\};

\draw[thick] (mtrx-2-1.north) -| (mtrx-5-1.south west)
                              -- (mtrx-5-1.south);
\draw[thick] (mtrx-2-4.north) -| (mtrx-5-4.south east)
                              -- (mtrx-5-4.south);
\path[draw,strip=orange]
    (mtrx-2-2.center) edge (mtrx-5-5.center)
    (mtrx-2-4.center) edge (mtrx-5-7.center)
    (mtrx-5-2.center) edge (mtrx-2-5.center)
    (mtrx-5-4.center) edge  (mtrx-2-7.center)
    (mtrx-2-8.center) edge (mtrx-5-11.center)
    (mtrx-2-10.center) edge (mtrx-5-13.center)
    (mtrx-5-8.center) edge (mtrx-2-11.center)
    (mtrx-5-10.center) edge  (mtrx-2-13.center)
    (mtrx-2-14.center) edge (mtrx-5-17.center)
    (mtrx-2-16.center) edge (mtrx-5-19.center)
    (mtrx-5-14.center) edge (mtrx-2-17.center)
    (mtrx-5-16.center)  --  (mtrx-2-19.center);
    
\path[draw,strip=cyan]
    (mtrx-2-1.center) edge (mtrx-5-4.center)
    (mtrx-2-3.center) edge (mtrx-5-6.center)
    (mtrx-5-1.center) edge (mtrx-2-4.center)
    (mtrx-5-3.center) edge (mtrx-2-6.center)
    (mtrx-2-7.center) edge (mtrx-5-10.center)
    (mtrx-2-9.center) edge (mtrx-5-12.center)
    (mtrx-5-7.center) edge (mtrx-2-10.center)
    (mtrx-5-9.center) edge (mtrx-2-12.center)
    (mtrx-2-13.center) edge (mtrx-5-16.center)
    (mtrx-2-15.center) edge (mtrx-5-18.center)
    (mtrx-5-13.center) edge (mtrx-2-16.center)
    (mtrx-5-15.center)  --  (mtrx-2-18.center);
\end{tikzpicture}
\end{center}
The sum of the products produced by the blue diagonals gives us $551$ and the sum of the products produced by the orange diagonals gives us $-411.$ Therefore the determinant is $551-411=140$.

We obtained the same result using the traditional method of expansions by minors. 
\end{example}

\section{The $5\times 5$ Determinant}\label{5x5}
   Consider the matrix
   \[M=\mat{a_{1}^1&a_{2}^1&a_{3}^1&a_{4}^1&a_{5}^1\\ 
            a_{1}^2&a_{2}^2&a_{3}^2&a_{4}^2&a_{5}^2\\
            a_{1}^3&a_{2}^3&a_{3}^3&a_{4}^3&a_{5}^3\\
            a_{1}^4&a_{2}^4&a_{3}^4&a_{4}^4&a_{5}^4\\
            a_{1}^5&a_{2}^5&a_{3}^5&a_{4}^5&a_{5}^5}.\]
            
    When computing the determinant of a $5\times5$ matrix, we obtain a sum of $5!=120$ products, 60 of them contribute to our positive sum and the other 60 contribute to the negative sum. Each of these products $p$ can be expressed as:
    \[
        p = a_{x}^1\cdot a_{y}^2\cdot a_{z}^3\cdot a_{u}^4\cdot a_{w}^5, \ \text{where}\ \{x,y,z,u,w\}=\{1,2,3,4,5\}.
    \] 
    Consider $\mathcal{P}$ be the set of all such products and $S_5$ the set of all permutations of $\{1,2,3,4,5\}$ see~\cite{Judson2020}.
    \begin{definition}
         For each product $p\in \mathcal{P},$ we define $\pi_p\in S_5$ by 
        % \[\pi_p(1)=x,\  \pi_p(2)=y,\ \pi_p(3)=z,\ \pi_p(4)=u,\ \pi_p(5)=w\] 
    % or more naturally using permutations notation 
    \[\pi(a_{x}^1\cdot a_{y}^2\cdot a_{z}^3\cdot a_{u}^4\cdot a_{w}^5)=\pi_p=\Perm{x&y&z&u&w}.\]
    \end{definition}
    Denote by $A_5$ (the alternating group) the even permutations in $S_5$. Using this notation, we can write (see~\cite{Lang1987LinearAlgebra}) the determinant of the matrix $M$ by 
    \[\det(M)=\sum_{p\in\mathcal{P}\,:\, \pi_p\in A_5} p-\sum_{p\in\mathcal{P}\,:\, \pi_p\notin A_5} p.\]
    In other words, the determinant of $M$ can be written as the sum of all the products $p\in\mathcal{P}$ such that $\pi_p$ is even minus the sum of all the products $p\in\mathcal{P}$ such that $\pi_p$ is odd.

    We could decompose $\sum_{p\in\mathcal{P}\,:\, \pi_p\in A_5} p$ as $P_1+P_2+P_3+P_4+P_5+P_6,$ where
    \begin{align*}
        P_1 & = \pro{1}{2}{3}{4}{5}+\pro{2}{3}{4}{5}{1}+\pro{3}{4}{5}{1}{2}+\pro{4}{5}{1}{2}{3}+\pro{5}{1}{2}{3}{4} \\
        & +\pro{5}{4}{3}{2}{1}+\pro{1}{5}{4}{3}{2}+\pro{2}{1}{5}{4}{3}+\pro{3}{2}{1}{5}{4}+\pro{4}{3}{2}{1}{5},\\
         P_2 & = \pro{4}{3}{5}{2}{1}+\pro{3}{5}{2}{1}{4}
         +\pro{5}{2}{1}{4}{3} +\pro{2}{1}{4}{3}{5}
         +\pro{1}{4}{3}{5}{2} \\
        & +\pro{1}{2}{5}{3}{4}+\pro{4}{1}{2}{5}{3}
        +\pro{3}{4}{1}{2}{5} +\pro{5}{3}{4}{1}{2}
          +\pro{2}{5}{3}{4}{1},\\
        P_3 & = \pro{2}{5}{1}{3}{4}+\pro{5}{1}{3}{4}{2}
     +\pro{1}{3}{4}{2}{5} +\pro{3}{4}{2}{5}{1}+ \pro{4}{2}{5}{1}{3}\\
    &  +\pro{4}{3}{1}{5}{2}+\pro{2}{4}{3}{1}{5}+\pro{5}{2}{4}{3}{1}
         +\pro{1}{5}{2}{4}{3}+\pro{3}{1}{5}{2}{4},\\
       P_4 & =  \pro{3}{1}{4}{5}{2}+\pro{1}{4}{5}{2}{3} 
         +\pro{4}{5}{2}{3}{1} +\pro{5}{2}{3}{1}{4} +\pro{2}{3}{1}{4}{5} \\
        & +\pro{2}{5}{4}{1}{3}+\pro{3}{2}{5}{4}{1} 
         +\pro{1}{3}{2}{5}{4} +\pro{4}{1}{3}{2}{5}+\pro{5}{4}{1}{3}{2},\\
     P_5 & = \pro{5}{4}{2}{1}{3}+\pro{4}{2}{1}{3}{5} 
      +\pro{2}{1}{3}{5}{4} +\pro{1}{3}{5}{4}{2} +\pro{3}{5}{4}{2}{1} \\
        & +\pro{3}{1}{2}{4}{5}+\pro{5}{3}{1}{2}{4} 
        +\pro{4}{5}{3}{1}{2} +\pro{2}{4}{5}{3}{1} +\pro{1}{2}{4}{5}{3},\\
     P_6 & = \pro{1}{4}{2}{3}{5}+\pro{4}{2}{3}{5}{1}
        +\pro{2}{3}{5}{1}{4} +\pro{3}{5}{1}{4}{2}+\pro{5}{1}{4}{2}{3}\\
        & +\pro{5}{3}{2}{4}{1}+\pro{1}{5}{3}{2}{4} 
        +\pro{4}{1}{5}{3}{2} +\pro{2}{4}{1}{5}{3}+\pro{3}{2}{4}{1}{5}.
    \end{align*}

    Observe that the first product in $P_1$ corresponds to $e,$ the second product corresponds to the cycle $c=(1\,2\,3\,4\,5),$ the next three correspond to $c^2,c^3,$ and $c^4.$ The remaining products in $P_1$ correspond to $(1\,5)(2\,4), (2\,5)(3\,4),(1\,2)(3\,5),(1\,3)(4\,5),(1\,4)(2\,3)$ respectively. All these are even permutations. Likewise the products in $P_2,..,P_6$ were found to correspond to even permutations.

    \medskip
    
    $P_1$ can be obtained by multiplying diagonally the terms found on the augmented matrix consisting of the original matrix $M$ and a copy of the first four rows placed on the right, as shown on the following matrix (we identify this augmented matrix with the identity permutation, as we are not modifying the original matrix $M$).

 \ \,\hspace{-0.5cm}\xxmini{0.5}{\begin{tikzpicture}[
strip/.style = {
    draw=#1,
    line width=1.25em, opacity=0.2,
    line cap=round ,% only if you like them ...
                            },
                    ]
\matrix (mtrx)  [matrix of math nodes,
                 column sep=0.75em,
                 nodes={text height=2ex,text width=2ex}]
{
a_1^1 & a_2^1 & a_3^1 & a_4^1 & a_5^1 & a_1^1 & a_2^1 & a_3^1 &
a_4^1  \\
a_1^2 & a_2^2 & a_3^2 & a_4^2 & a_5^2 & a_1^2 & a_2^2 & a_3^2 &
a_4^2   \\
a_1^3 & a_2^3 & a_3^3 & a_4^3 & a_5^3 & a_1^3 & a_2^3 & a_3^3 &
a_4^3   \\
a_1^4 & a_2^4 & a_3^4 & a_4^4 & a_5^4 & a_1^4 & a_2^4 & a_3^4 &
a_4^4   \\
a_1^5 & a_2^5 & a_3^5 & a_4^5 & a_5^5 & a_1^5 & a_2^5 & a_3^5 &
a_4^5   \\
};
\draw[thick] (mtrx-1-1.north) -| (mtrx-5-1.south west)
                              -- (mtrx-5-1.south);
\draw[thick] (mtrx-1-5.north) -| (mtrx-5-5.south east)
                              -- (mtrx-5-5.south);
\draw [decorate,decoration={brace,amplitude=10pt,raise=2pt},yshift=0pt]
(mtrx-1-1.north west) -- (mtrx-1-5.north east) node [black,midway,yshift=.6cm] {\footnotesize
corresponds to  $e=$\perm{1&2&3&4&5}};
\draw [decorate,decoration={brace,mirror,amplitude=10pt},xshift=3pt]
(mtrx-5-6.south west) -- (mtrx-5-9.south east) node [black,midway,yshift=-.5cm] {\footnotesize
copy first four columns};
\path[draw,strip=cyan,shorten <=-1.5mm,shorten >=-1.5mm]
    (mtrx-1-1.center) edge (mtrx-5-5.center)
    (mtrx-1-2.center) edge (mtrx-5-6.center)
    (mtrx-1-3.center) edge (mtrx-5-7.center)
    (mtrx-1-4.center) edge (mtrx-5-8.center)
    (mtrx-1-5.center) edge (mtrx-5-9.center); %
\path[draw,strip=cyan,shorten <=-1.5mm,shorten >=-1.5mm]
    (mtrx-5-1.center) edge (mtrx-1-5.center)
    (mtrx-5-2.center) edge (mtrx-1-6.center)
    (mtrx-5-3.center) edge (mtrx-1-7.center)
    (mtrx-5-4.center) edge (mtrx-1-8.center)
    (mtrx-5-5.center) edge (mtrx-1-9.center);
\end{tikzpicture}}\hspace{0cm}
   \xxmini{0.5}{\small \begin{align*}
        P_1 & = \pro{1}{2}{3}{4}{5}+\pro{2}{3}{4}{5}{1}+\pro{3}{4}{5}{1}{2}\\
        & +\pro{4}{5}{1}{2}{3}+\pro{5}{1}{2}{3}{4} \\
        & +\pro{5}{4}{3}{2}{1}+\pro{1}{5}{4}{3}{2}+\pro{2}{1}{5}{4}{3}\\
        & +\pro{3}{2}{1}{5}{4}+\pro{4}{3}{2}{1}{5}.\\
    \end{align*}}

    \bigskip
    
    Likewise, $P_2$ can be obtained by multiplying diagonally in the following augmented matrix (which we identify with \perm{4&3&5&2&1} as it corresponds to rearranging the columns of $M$ in the following order, column 4, column 3, column 5, column 2 and column 1).

   \ \,\hspace{-0.5cm}\xxmini{0.5}{\begin{tikzpicture}[
strip/.style = {
    draw=#1,
    line width=1.25em, opacity=0.2,
    line cap=round ,% only if you like them ...
                            },
                    ]
\matrix (mtrx)  [matrix of math nodes,
                 column sep=0.75em,
                 nodes={text height=2ex,text width=2ex}]
{
a_4^1 & a_3^1 & a_5^1 & a_2^1 & a_1^1 & a_4^1 & a_3^1 & a_5^1 & a_2^1  \\
a_4^2 & a_3^2 & a_5^2 & a_2^2 & a_1^2 & a_4^2 & a_3^2 & a_5^2 & a_2^2  \\
a_4^3 & a_3^3 & a_5^3 & a_2^3 & a_1^3 & a_4^3 & a_3^3 & a_5^3 & a_2^3  \\
a_4^4 & a_3^4 & a_5^4 & a_2^4 & a_1^4 & a_4^4 & a_3^4 & a_5^4 & a_2^4  \\
a_4^5 & a_3^5 & a_5^5 & a_2^5 & a_1^5 & a_4^5 & a_3^5 & a_5^5 & a_2^5  \\
};
\draw[thick] (mtrx-1-1.north) -| (mtrx-5-1.south west)
                              -- (mtrx-5-1.south);
\draw[thick] (mtrx-1-5.north) -| (mtrx-5-5.south east)
                              -- (mtrx-5-5.south);
\draw [decorate,decoration={brace,amplitude=10pt,raise=2pt},yshift=0pt]
(mtrx-1-1.north west) -- (mtrx-1-5.north east) node [black,midway,yshift=.6cm] {\footnotesize
corresponds to  $\sigma=$\perm{4&3&5&2&1}};
\draw [decorate,decoration={brace,mirror,amplitude=10pt},xshift=3pt]
(mtrx-5-6.south west) -- (mtrx-5-9.south east) node [black,midway,yshift=-.5cm] {\footnotesize
copy first four columns};
\path[draw,strip=cyan,shorten <=-1.5mm,shorten >=-1.5mm]
    (mtrx-1-1.center) edge (mtrx-5-5.center)
    (mtrx-1-2.center) edge (mtrx-5-6.center)
    (mtrx-1-3.center) edge (mtrx-5-7.center)
    (mtrx-1-4.center) edge (mtrx-5-8.center)
    (mtrx-1-5.center) edge (mtrx-5-9.center); %
\path[draw,strip=cyan,shorten <=-1.5mm,shorten >=-1.5mm]
    (mtrx-5-1.center) edge (mtrx-1-5.center)
    (mtrx-5-2.center) edge (mtrx-1-6.center)
    (mtrx-5-3.center) edge (mtrx-1-7.center)
    (mtrx-5-4.center) edge (mtrx-1-8.center)
    (mtrx-5-5.center) edge (mtrx-1-9.center);
\end{tikzpicture}}\hspace{0cm}
   \xxmini{0.5}{\small \begin{align*}
        P_2 & = \pro{4}{3}{5}{2}{1}+\pro{3}{5}{2}{1}{4}
         +\pro{5}{2}{1}{4}{3}\\
        & +\pro{2}{1}{4}{3}{5}
         +\pro{1}{4}{3}{5}{2} \\
        & +\pro{1}{2}{5}{3}{4}+\pro{4}{1}{2}{5}{3}
        +\pro{3}{4}{1}{2}{5}\\
        & +\pro{5}{3}{4}{1}{2}
          +\pro{2}{5}{3}{4}{1}
    \end{align*}}

    Similarly $P_3,P_4,P_5$ and $P_6$ can be obtained by multiplying diagonally in the following augmented matrices (we identify each one with a specific permutation corresponding to the first five columns). 

    \ \,\hspace{-0.5cm}\xxmini{0.5}{\begin{tikzpicture}[
strip/.style = {
    draw=#1,
    line width=1.25em, opacity=0.2,
    line cap=round ,% only if you like them ...
                            },
                    ]
\matrix (mtrx)  [matrix of math nodes,
                 column sep=0.75em,
                 nodes={text height=2ex,text width=2ex}]
{
a_2^1 & a_5^1 & a_1^1 & a_3^1 & a_4^1 & a_2^1 & a_5^1 & a_1^1 & a_3^1  \\
a_2^2 & a_5^2 & a_1^2 & a_3^2 & a_4^2 & a_2^2 & a_5^2 & a_1^2 & a_3^2  \\
a_2^3 & a_5^3 & a_1^3 & a_3^3 & a_4^3 & a_2^3 & a_5^3 & a_1^3 & a_3^3  \\
a_2^4 & a_5^4 & a_1^4 & a_3^4 & a_4^4 & a_2^4 & a_5^4 & a_1^4 & a_3^4  \\
a_2^5 & a_5^5 & a_1^5 & a_3^5 & a_4^5 & a_2^5 & a_5^5 & a_1^5 & a_3^5  \\
};
\draw[thick] (mtrx-1-1.north) -| (mtrx-5-1.south west)
                              -- (mtrx-5-1.south);
\draw[thick] (mtrx-1-5.north) -| (mtrx-5-5.south east)
                              -- (mtrx-5-5.south);
\draw [decorate,decoration={brace,amplitude=10pt,raise=2pt},yshift=0pt]
(mtrx-1-1.north west) -- (mtrx-1-5.north east) node [black,midway,yshift=.6cm] {\footnotesize
corresponds to  $\sigma^2=$\perm{2&5&1&3&4}};
\draw [decorate,decoration={brace,mirror,amplitude=10pt},xshift=3pt]
(mtrx-5-6.south west) -- (mtrx-5-9.south east) node [black,midway,yshift=-.5cm] {\footnotesize
copy first four columns};
\path[draw,strip=cyan,shorten <=-1.5mm,shorten >=-1.5mm]
    (mtrx-1-1.center) edge (mtrx-5-5.center)
    (mtrx-1-2.center) edge (mtrx-5-6.center)
    (mtrx-1-3.center) edge (mtrx-5-7.center)
    (mtrx-1-4.center) edge (mtrx-5-8.center)
    (mtrx-1-5.center) edge (mtrx-5-9.center); %
\path[draw,strip=cyan,shorten <=-1.5mm,shorten >=-1.5mm]
    (mtrx-5-1.center) edge (mtrx-1-5.center)
    (mtrx-5-2.center) edge (mtrx-1-6.center)
    (mtrx-5-3.center) edge (mtrx-1-7.center)
    (mtrx-5-4.center) edge (mtrx-1-8.center)
    (mtrx-5-5.center) edge (mtrx-1-9.center);
\end{tikzpicture}}\hspace{0cm}
   \xxmini{0.5}{\small \begin{align*}
    P_3 & = \pro{2}{5}{1}{3}{4}+\pro{5}{1}{3}{4}{2}
     +\pro{1}{3}{4}{2}{5}\\
    & +\pro{3}{4}{2}{5}{1}+ \pro{4}{2}{5}{1}{3}\\
    &  +\pro{4}{3}{1}{5}{2}+\pro{2}{4}{3}{1}{5}+\pro{5}{2}{4}{3}{1}\\
        &  +\pro{1}{5}{2}{4}{3}+\pro{3}{1}{5}{2}{4}.
    \end{align*}}
    
  \bigskip

    \ \,\hspace{-0.5cm}\xxmini{0.5}{\begin{tikzpicture}[
strip/.style = {
    draw=#1,
    line width=1.25em, opacity=0.2,
    line cap=round ,% only if you like them ...
                            },
                    ]
\matrix (mtrx)  [matrix of math nodes,
                 column sep=0.75em,
                 nodes={text height=2ex,text width=2ex}]
{
a_3^1 & a_1^1 & a_4^1 & a_5^1 & a_2^1 & a_3^1 & a_1^1 & a_4^1 & a_5^1  \\
a_3^2 & a_1^2 & a_4^2 & a_5^2 & a_2^2 & a_3^2 & a_1^2 & a_4^2 & a_5^2  \\
a_3^3 & a_1^3 & a_4^3 & a_5^3 & a_2^3 & a_3^3 & a_1^3 & a_4^3 & a_5^3  \\
a_3^4 & a_1^4 & a_4^4 & a_5^4 & a_2^4 & a_3^4 & a_1^4 & a_4^4 & a_5^4  \\
a_3^5 & a_1^5 & a_4^5 & a_5^5 & a_2^5 & a_3^5 & a_1^5 & a_4^5 & a_5^5  \\
};
\draw[thick] (mtrx-1-1.north) -| (mtrx-5-1.south west)
                              -- (mtrx-5-1.south);
\draw[thick] (mtrx-1-5.north) -| (mtrx-5-5.south east)
                              -- (mtrx-5-5.south);
\draw [decorate,decoration={brace,amplitude=10pt,raise=2pt},yshift=0pt]
(mtrx-1-1.north west) -- (mtrx-1-5.north east) node [black,midway,yshift=.6cm] {\footnotesize
corresponds to  $\sigma^3=$\perm{3&1&4&5&2}};
\draw [decorate,decoration={brace,mirror,amplitude=10pt},xshift=3pt]
(mtrx-5-6.south west) -- (mtrx-5-9.south east) node [black,midway,yshift=-.5cm] {\footnotesize
copy first four columns};
\path[draw,strip=cyan,shorten <=-1.5mm,shorten >=-1.5mm]
    (mtrx-1-1.center) edge (mtrx-5-5.center)
    (mtrx-1-2.center) edge (mtrx-5-6.center)
    (mtrx-1-3.center) edge (mtrx-5-7.center)
    (mtrx-1-4.center) edge (mtrx-5-8.center)
    (mtrx-1-5.center) edge (mtrx-5-9.center); %
\path[draw,strip=cyan,shorten <=-1.5mm,shorten >=-1.5mm]
    (mtrx-5-1.center) edge (mtrx-1-5.center)
    (mtrx-5-2.center) edge (mtrx-1-6.center)
    (mtrx-5-3.center) edge (mtrx-1-7.center)
    (mtrx-5-4.center) edge (mtrx-1-8.center)
    (mtrx-5-5.center) edge (mtrx-1-9.center);
\end{tikzpicture}}\hspace{0cm}
   \xxmini{0.5}{\small \begin{align*}
        P_4 & =  \pro{3}{1}{4}{5}{2}+\pro{1}{4}{5}{2}{3} 
         +\pro{4}{5}{2}{3}{1}\\
        & +\pro{5}{2}{3}{1}{4} +\pro{2}{3}{1}{4}{5} \\
        & +\pro{2}{5}{4}{1}{3}+\pro{3}{2}{5}{4}{1} 
         +\pro{1}{3}{2}{5}{4}\\
        & +\pro{4}{1}{3}{2}{5}+\pro{5}{4}{1}{3}{2}.
    \end{align*}}
    
   \bigskip

       \ \,\hspace{-0.5cm}\xxmini{0.5}{\begin{tikzpicture}[
strip/.style = {
    draw=#1,
    line width=1.25em, opacity=0.2,
    line cap=round ,% only if you like them ...
                            },
                    ]
\matrix (mtrx)  [matrix of math nodes,
                 column sep=0.75em,
                 nodes={text height=2ex,text width=2ex}]
{
a_5^1 & a_4^1 & a_2^1 & a_1^1 & a_3^1 & a_5^1 & a_4^1 & a_2^1 & a_1^1  \\
a_5^2 & a_4^2 & a_2^2 & a_1^2 & a_3^2 & a_5^2 & a_4^2 & a_2^2 & a_1^2  \\
a_5^3 & a_4^3 & a_2^3 & a_1^3 & a_3^3 & a_5^3 & a_4^3 & a_2^3 & a_1^3  \\
a_5^4 & a_4^4 & a_2^4 & a_1^4 & a_3^4 & a_5^4 & a_4^4 & a_2^4 & a_1^4  \\
a_5^5 & a_4^5 & a_2^5 & a_1^5 & a_3^5 & a_5^5 & a_4^5 & a_2^5 & a_1^5  \\
};
\draw[thick] (mtrx-1-1.north) -| (mtrx-5-1.south west)
                              -- (mtrx-5-1.south);
\draw[thick] (mtrx-1-5.north) -| (mtrx-5-5.south east)
                              -- (mtrx-5-5.south);
\draw [decorate,decoration={brace,amplitude=10pt,raise=2pt},yshift=0pt]
(mtrx-1-1.north west) -- (mtrx-1-5.north east) node [black,midway,yshift=.6cm] {\footnotesize
corresponds to  $\sigma^4=$\perm{5&4&2&1&3}};
\draw [decorate,decoration={brace,mirror,amplitude=10pt},xshift=3pt]
(mtrx-5-6.south west) -- (mtrx-5-9.south east) node [black,midway,yshift=-.5cm] {\footnotesize
copy first four columns};
\path[draw,strip=cyan,shorten <=-1.5mm,shorten >=-1.5mm]
    (mtrx-1-1.center) edge (mtrx-5-5.center)
    (mtrx-1-2.center) edge (mtrx-5-6.center)
    (mtrx-1-3.center) edge (mtrx-5-7.center)
    (mtrx-1-4.center) edge (mtrx-5-8.center)
    (mtrx-1-5.center) edge (mtrx-5-9.center); %
\path[draw,strip=cyan,shorten <=-1.5mm,shorten >=-1.5mm]
    (mtrx-5-1.center) edge (mtrx-1-5.center)
    (mtrx-5-2.center) edge (mtrx-1-6.center)
    (mtrx-5-3.center) edge (mtrx-1-7.center)
    (mtrx-5-4.center) edge (mtrx-1-8.center)
    (mtrx-5-5.center) edge (mtrx-1-9.center);
\end{tikzpicture}}\hspace{0cm}
   \xxmini{0.5}{\small \begin{align*}
    P_5 & = \pro{5}{4}{2}{1}{3}+\pro{4}{2}{1}{3}{5} 
      +\pro{2}{1}{3}{5}{4}\\
        & +\pro{1}{3}{5}{4}{2} +\pro{3}{5}{4}{2}{1} \\
        & +\pro{3}{1}{2}{4}{5}+\pro{5}{3}{1}{2}{4} 
        +\pro{4}{5}{3}{1}{2}\\
        & +\pro{2}{4}{5}{3}{1} +\pro{1}{2}{4}{5}{3}.
    \end{align*}}
    
   \bigskip

       \ \,\hspace{-0.5cm}\xxmini{0.5}{\begin{tikzpicture}[
strip/.style = {
    draw=#1,
    line width=1.25em, opacity=0.2,
    line cap=round ,% only if you like them ...
                            },
                    ]
\matrix (mtrx)  [matrix of math nodes,
                 column sep=0.75em,
                 nodes={text height=2ex,text width=2ex}]
{
a_1^1 & a_4^1 & a_2^1 & a_3^1 & a_5^1 & a_1^1 & a_4^1 & a_2^1 & a_3^1  \\
a_1^2 & a_4^2 & a_2^2 & a_3^2 & a_5^2 & a_1^2 & a_4^2 & a_2^2 & a_3^2  \\
a_1^3 & a_4^3 & a_2^3 & a_3^3 & a_5^3 & a_1^3 & a_4^3 & a_2^3 & a_3^3  \\
a_1^4 & a_4^4 & a_2^4 & a_3^4 & a_5^4 & a_1^4 & a_4^4 & a_2^4 & a_3^4  \\
a_1^5 & a_4^5 & a_2^5 & a_3^5 & a_5^5 & a_1^5 & a_4^5 & a_2^5 & a_3^5  \\
};
\draw[thick] (mtrx-1-1.north) -| (mtrx-5-1.south west)
                              -- (mtrx-5-1.south);
\draw[thick] (mtrx-1-5.north) -| (mtrx-5-5.south east)
                              -- (mtrx-5-5.south);
\draw [decorate,decoration={brace,amplitude=10pt,raise=2pt},yshift=0pt]
(mtrx-1-1.north west) -- (mtrx-1-5.north east) node [black,midway,yshift=.6cm] {\footnotesize
corresponds to  $\tau=$\perm{1&4&2&3&5}};
\draw [decorate,decoration={brace,mirror,amplitude=10pt},xshift=3pt]
(mtrx-5-6.south west) -- (mtrx-5-9.south east) node [black,midway,yshift=-.5cm] {\footnotesize
copy first four columns};
\path[draw,strip=cyan,shorten <=-1.5mm,shorten >=-1.5mm]
    (mtrx-1-1.center) edge (mtrx-5-5.center)
    (mtrx-1-2.center) edge (mtrx-5-6.center)
    (mtrx-1-3.center) edge (mtrx-5-7.center)
    (mtrx-1-4.center) edge (mtrx-5-8.center)
    (mtrx-1-5.center) edge (mtrx-5-9.center); %
\path[draw,strip=cyan,shorten <=-1.5mm,shorten >=-1.5mm]
    (mtrx-5-1.center) edge (mtrx-1-5.center)
    (mtrx-5-2.center) edge (mtrx-1-6.center)
    (mtrx-5-3.center) edge (mtrx-1-7.center)
    (mtrx-5-4.center) edge (mtrx-1-8.center)
    (mtrx-5-5.center) edge (mtrx-1-9.center);
\end{tikzpicture}}\hspace{0cm}
   \xxmini{0.5}{\small \begin{align*}
        P_6 & = \pro{1}{4}{2}{3}{5}+\pro{4}{2}{3}{5}{1}
        +\pro{2}{3}{5}{1}{4}\\
        & +\pro{3}{5}{1}{4}{2}+\pro{5}{1}{4}{2}{3}\\
        & +\pro{5}{3}{2}{4}{1}+\pro{1}{5}{3}{2}{4} 
        +\pro{4}{1}{5}{3}{2}\\
        & +\pro{2}{4}{1}{5}{3}+\pro{3}{2}{4}{1}{5}.
    \end{align*}}

When putting together all these augmented matrices--which give the products corresponding to all the even permutations--they connect in such a way that the last column of one is the first column of the next one.  By placing the first two matrices together, we observe that they connected similarly to how we connected the $4\times 4$ matrix.

 \ \,\hspace{-0.5cm}\xxmini{0.5}{\begin{tikzpicture}[
strip/.style = {
    draw=#1,
    line width=1.25em, opacity=0.2,
    line cap=round ,% only if you like them ...
                            },
                    ]
\matrix (mtrx)  [matrix of math nodes,
                 column sep=1em,
                 nodes={text height=2ex,text width=2ex}]
{
a_1^1 & a_2^1 & a_3^1 & a_4^1 & a_5^1 & a_1^1 & a_2^1 & a_3^1 &
a_4^1 & a_3^1 & a_5^1 & a_2^1 & a_1^1 & a_4^1 & a_3^1 & a_5^1 & a_2^1 \\
a_1^2 & a_2^2 & a_3^2 & a_4^2 & a_5^2 & a_1^2 & a_2^2 & a_3^2 &
a_4^2 & a_3^2 & a_5^2 & a_2^2 & a_1^2 & a_4^2 & a_3^2 & a_5^2 & a_2^2   \\
a_1^3 & a_2^3 & a_3^3 & a_4^3 & a_5^3 & a_1^3 & a_2^3 & a_3^3 &
a_4^3 & a_3^3 & a_5^3 & a_2^3 & a_1^3 & a_4^3 & a_3^3 & a_5^3 & a_2^3   \\
a_1^4 & a_2^4 & a_3^4 & a_4^4 & a_5^4 & a_1^4 & a_2^4 & a_3^4 &
a_4^4 & a_3^4 & a_5^4 & a_2^4 & a_1^4 & a_4^4 & a_3^4 & a_5^4 & a_2^4   \\
a_1^5 & a_2^5 & a_3^5 & a_4^5 & a_5^5 & a_1^5 & a_2^5 & a_3^5 &
a_4^5 & a_3^5 & a_5^5 & a_2^5 & a_1^5 & a_4^5 & a_3^5 & a_5^5 & a_2^5   \\
};
\draw[thick] (mtrx-1-1.north) -| (mtrx-5-1.south west)
                              -- (mtrx-5-1.south);
\draw[thick] (mtrx-1-5.north) -| (mtrx-5-5.south east)
                              -- (mtrx-5-5.south);
\draw [decorate,decoration={brace,amplitude=10pt,raise=2pt},yshift=0pt]
(mtrx-1-1.north west) -- (mtrx-1-5.north east) node [black,midway,yshift=.6cm] {\footnotesize
corresponds to  $e=$\perm{1&2&3&4&5}}; 
\path[draw,strip=cyan,shorten <=-1.5mm,shorten >=-1.5mm]
    (mtrx-1-1.center) edge (mtrx-5-5.center)
    (mtrx-1-2.center) edge (mtrx-5-6.center)
    (mtrx-1-3.center) edge (mtrx-5-7.center)
    (mtrx-1-4.center) edge (mtrx-5-8.center)
    (mtrx-1-5.center) edge (mtrx-5-9.center); %
\path[draw,strip=cyan,shorten <=-1.5mm,shorten >=-1.5mm]
    (mtrx-5-1.center) edge (mtrx-1-5.center)
    (mtrx-5-2.center) edge (mtrx-1-6.center)
    (mtrx-5-3.center) edge (mtrx-1-7.center)
    (mtrx-5-4.center) edge (mtrx-1-8.center)
    (mtrx-5-5.center) edge (mtrx-1-9.center);
\draw [decorate,decoration={brace,amplitude=10pt,raise=2pt},yshift=0pt]
(mtrx-1-9.north west) -- (mtrx-1-13.north east) node [black,midway,yshift=.6cm] {\footnotesize
corresponds to  $\sigma=$\perm{4&3&5&2&1}}; 
\path[draw,strip=cyan,shorten <=-1.5mm,shorten >=-1.5mm]
    (mtrx-1-9.center) edge (mtrx-5-13.center)
    (mtrx-1-10.center) edge (mtrx-5-14.center)
    (mtrx-1-11.center) edge (mtrx-5-15.center)
    (mtrx-1-12.center) edge (mtrx-5-16.center)
    (mtrx-1-13.center) edge (mtrx-5-17.center); %
\path[draw,strip=cyan,shorten <=-1.5mm,shorten >=-1.5mm]
    (mtrx-5-9.center) edge (mtrx-1-13.center)
    (mtrx-5-10.center) edge (mtrx-1-14.center)
    (mtrx-5-11.center) edge (mtrx-1-15.center)
    (mtrx-5-12.center) edge (mtrx-1-16.center)
    (mtrx-5-13.center) edge (mtrx-1-17.center);
\end{tikzpicture}}

\bigskip

We can further connect all of the them together creating a much larger augmented matrix that looks like a quilt:

\medskip
 
\ \,\hspace{-0.65cm}\includegraphics[width=0.75\paperwidth]{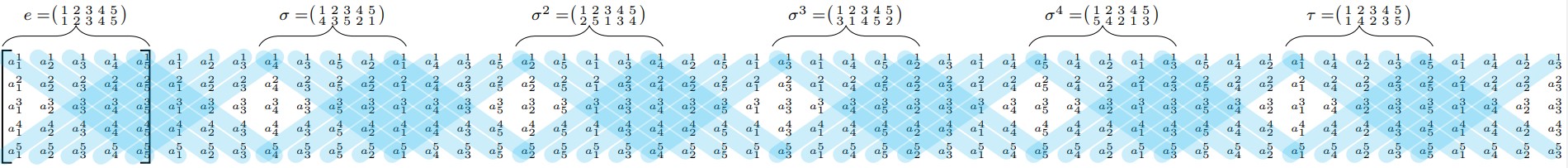}

When looking at the products corresponding to odd permutations, a very similar process occurs; here we decompose $\sum_{p\in\mathcal{P}\,:\, \pi_p\notin A_5} p$ as $N_1+N_2+N_3+N_4+N_5+N_6,$ where
    \begin{align*}
        N_1 & = \pro{1}{2}{3}{4}{5}+\pro{2}{3}{4}{5}{1}+\pro{3}{4}{5}{1}{2}+\pro{4}{5}{1}{2}{3}+\pro{5}{1}{2}{3}{4} \\
        & +\pro{5}{4}{3}{2}{1}+\pro{1}{5}{4}{3}{2}+\pro{2}{1}{5}{4}{3}+\pro{3}{2}{1}{5}{4}+\pro{4}{3}{2}{1}{5},\\
        N_2 & =  \pro{3}{4}{5}{2}{1}+\pro{4}{5}{2}{1}{3} +\pro{5}{2}{1}{3}{4}+\pro{2}{1}{3}{4}{5} +\pro{1}{3}{4}{5}{2} \\
        & +\pro{1}{2}{5}{4}{3}+\pro{3}{1}{2}{5}{4} +\pro{4}{3}{1}{2}{5}+\pro{5}{4}{3}{1}{2}+\pro{2}{5}{4}{3}{1},\\
        N_3 & =  \pro{2}{5}{1}{4}{3}+\pro{5}{1}{4}{3}{2}
        + \pro{1}{4}{3}{2}{5}+\pro{4}{3}{2}{5}{1}
        + \pro{3}{2}{5}{1}{4} \\
        & +\pro{3}{4}{1}{5}{2}+\pro{2}{3}{4}{1}{5}
        +\pro{5}{2}{3}{4}{1}+\pro{1}{5}{2}{3}{4}
        +\pro{4}{1}{5}{2}{3},\\
        N_4 & = \pro{4}{1}{3}{5}{2}+\pro{1}{3}{5}{2}{4} 
        +\pro{3}{5}{2}{4}{1}+\pro{5}{2}{4}{1}{3} 
        +\pro{2}{4}{1}{3}{5} \\
        & +\pro{2}{5}{3}{1}{4}+\pro{4}{2}{5}{3}{1} 
          +\pro{1}{4}{2}{5}{3}+\pro{3}{1}{4}{2}{5} 
          +\pro{5}{3}{1}{4}{2}, \\
        N_5 & =  \pro{5}{3}{2}{1}{4}+\pro{3}{2}{1}{4}{5} 
        +\pro{2}{1}{4}{5}{3}+\pro{1}{4}{5}{3}{2} 
        +\pro{4}{5}{3}{2}{1} \\
        & +\pro{4}{1}{2}{3}{5}+\pro{5}{4}{1}{2}{3} 
         +\pro{3}{5}{4}{1}{2}+\pro{2}{3}{5}{4}{1} 
         +\pro{1}{2}{3}{5}{4},\\
        N_6 & =   \pro{1}{3}{2}{4}{5}+\pro{3}{2}{4}{5}{1}
         +\pro{2}{4}{5}{1}{3}+\pro{4}{5}{1}{3}{2}
         +\pro{5}{1}{3}{2}{4}\\
        & +\pro{5}{4}{2}{3}{1}+\pro{1}{5}{4}{2}{3}
         +\pro{3}{1}{5}{4}{2}+\pro{2}{3}{1}{5}{4}
         +\pro{4}{2}{3}{1}{5}.
    \end{align*}

The products that appear in $N_1,...,N_6$ are the same products that appear in $P_1,...,P_6$ multiplied by $(3\,4)$ which result in odd permutations.
These products are obtained by multiplying diagonally the following matrices:

 \ \,\hspace{-0.75cm}\xxmini{0.5}{\begin{tikzpicture}[
strip/.style = {
    draw=#1,
    line width=1.25em, opacity=0.2,
    line cap=round ,% only if you like them ...
                            },
                    ]
\matrix (mtrx)  [matrix of math nodes,
                 column sep=0.75em,
                 nodes={text height=2ex,text width=2ex}]
{
a_1^1 & a_2^1 & a_4^1 & a_3^1 & a_5^1 & a_1^1 & a_2^1 & a_4^1 &
a_3^1  \\
a_1^2 & a_2^2 & a_4^2 & a_3^2 & a_5^2 & a_1^2 & a_2^2 & a_4^2 &
a_3^2   \\
a_1^3 & a_2^3 & a_4^3 & a_3^3 & a_5^3 & a_1^3 & a_2^3 & a_4^3 &
a_3^3   \\
a_1^4 & a_2^4 & a_4^4 & a_3^4 & a_5^4 & a_1^4 & a_2^4 & a_4^4 &
a_3^4   \\
a_1^5 & a_2^5 & a_4^5 & a_3^5 & a_5^5 & a_1^5 & a_2^5 & a_4^5 &
a_3^5   \\ 
};
\draw[thick] (mtrx-1-1.north) -| (mtrx-5-1.south west)
                              -- (mtrx-5-1.south);
\draw[thick] (mtrx-1-5.north) -| (mtrx-5-5.south east)
                              -- (mtrx-5-5.south);
\draw [decorate,decoration={brace,amplitude=10pt,raise=2pt},yshift=0pt]
(mtrx-1-1.north west) -- (mtrx-1-5.north east) node [black,midway,yshift=.6cm] {\footnotesize
corresponds to  $(34)e=$\perm{1&2&4&3&5}};
\draw [decorate,decoration={brace,mirror,amplitude=10pt},xshift=3pt]
(mtrx-5-6.south west) -- (mtrx-5-9.south east) node [black,midway,yshift=-.5cm] {\footnotesize
copy first four columns};
\path[draw,strip=orange,shorten <=-1.5mm,shorten >=-1.5mm]
    (mtrx-1-1.center) edge (mtrx-5-5.center)
    (mtrx-1-2.center) edge (mtrx-5-6.center)
    (mtrx-1-3.center) edge (mtrx-5-7.center)
    (mtrx-1-4.center) edge (mtrx-5-8.center)
    (mtrx-1-5.center) edge (mtrx-5-9.center); %
\path[draw,strip=orange,shorten <=-1.5mm,shorten >=-1.5mm]
    (mtrx-5-1.center) edge (mtrx-1-5.center)
    (mtrx-5-2.center) edge (mtrx-1-6.center)
    (mtrx-5-3.center) edge (mtrx-1-7.center)
    (mtrx-5-4.center) edge (mtrx-1-8.center)
    (mtrx-5-5.center) edge (mtrx-1-9.center);
\end{tikzpicture}}\hspace{0cm}
   \xxmini{0.5}{\small \begin{align*}
        N_1 & = \pro{1}{2}{4}{3}{5}+\pro{2}{4}{3}{5}{1}
         +\pro{4}{3}{5}{1}{2}\\
        & +\pro{3}{5}{1}{2}{4} +\pro{5}{1}{2}{4}{3} \\
        & +\pro{5}{3}{4}{2}{1}+\pro{1}{5}{3}{4}{2}
        +\pro{2}{1}{5}{3}{4}\\
        & +\pro{4}{2}{1}{5}{3}+\pro{3}{4}{2}{1}{5}.
    \end{align*}}

 \ \,\hspace{-0.75cm}\xxmini{0.5}{\begin{tikzpicture}[
strip/.style = {
    draw=#1,
    line width=1.25em, opacity=0.2,
    line cap=round ,% only if you like them ...
                            },
                    ]
\matrix (mtrx)  [matrix of math nodes,
                 column sep=0.75em,
                 nodes={text height=2ex,text width=2ex}]
{
a_3^1 & a_4^1 & a_5^1 & a_2^1 & a_1^1 & a_3^1 & a_4^1 & a_5^1 & a_2^1  \\
a_3^2 & a_4^2 & a_5^2 & a_2^2 & a_1^2 & a_3^2 & a_4^2 & a_5^2 & a_2^2  \\
a_3^3 & a_4^3 & a_5^3 & a_2^3 & a_1^3 & a_3^3 & a_4^3 & a_5^3 & a_2^3  \\
a_3^4 & a_4^4 & a_5^4 & a_2^4 & a_1^4 & a_3^4 & a_4^4 & a_5^4 & a_2^4  \\
a_3^5 & a_4^5 & a_5^5 & a_2^5 & a_1^5 & a_3^5 & a_4^5 & a_5^5 & a_2^5  \\
};
\draw[thick] (mtrx-1-1.north) -| (mtrx-5-1.south west)
                              -- (mtrx-5-1.south);
\draw[thick] (mtrx-1-5.north) -| (mtrx-5-5.south east)
                              -- (mtrx-5-5.south);
\draw [decorate,decoration={brace,amplitude=10pt,raise=2pt},yshift=0pt]
(mtrx-1-1.north west) -- (mtrx-1-5.north east) node [black,midway,yshift=.6cm] {\footnotesize
corresponds to  $(34)\cdot\sigma=$\perm{3&4&5&2&1}};
\draw [decorate,decoration={brace,mirror,amplitude=10pt},xshift=3pt]
(mtrx-5-6.south west) -- (mtrx-5-9.south east) node [black,midway,yshift=-.5cm] {\footnotesize
copy first four columns};
\path[draw,strip=orange,shorten <=-1.5mm,shorten >=-1.5mm]
    (mtrx-1-1.center) edge (mtrx-5-5.center)
    (mtrx-1-2.center) edge (mtrx-5-6.center)
    (mtrx-1-3.center) edge (mtrx-5-7.center)
    (mtrx-1-4.center) edge (mtrx-5-8.center)
    (mtrx-1-5.center) edge (mtrx-5-9.center); %
\path[draw,strip=orange,shorten <=-1.5mm,shorten >=-1.5mm]
    (mtrx-5-1.center) edge (mtrx-1-5.center)
    (mtrx-5-2.center) edge (mtrx-1-6.center)
    (mtrx-5-3.center) edge (mtrx-1-7.center)
    (mtrx-5-4.center) edge (mtrx-1-8.center)
    (mtrx-5-5.center) edge (mtrx-1-9.center);
\end{tikzpicture}}\hspace{0cm}
   \xxmini{0.5}{\small \begin{align*}
        N_2 & = \pro{3}{4}{5}{2}{1}+\pro{4}{5}{2}{1}{3}
         +\pro{5}{2}{1}{3}{4}\\
        & +\pro{2}{1}{3}{4}{5}+\pro{1}{3}{4}{5}{2} \\
        & +\pro{1}{2}{5}{4}{3}+\pro{3}{1}{2}{5}{4} 
        +\pro{4}{3}{1}{2}{5}\\
        & +\pro{5}{4}{3}{1}{2} +\pro{2}{5}{4}{3}{1}.
    \end{align*}}

    \ \,\hspace{-0.75cm}\xxmini{0.5}{\begin{tikzpicture}[
strip/.style = {
    draw=#1,
    line width=1.25em, opacity=0.2,
    line cap=round ,% only if you like them ...
                            },
                    ]
\matrix (mtrx)  [matrix of math nodes,
                 column sep=0.75em,
                 nodes={text height=2ex,text width=2ex}]
{
a_2^1 & a_5^1 & a_1^1 & a_4^1 & a_3^1 & a_2^1 & a_5^1 & a_1^1 & a_4^1  \\
a_2^2 & a_5^2 & a_1^2 & a_4^2 & a_3^2 & a_2^2 & a_5^2 & a_1^2 & a_4^2  \\
a_2^3 & a_5^3 & a_1^3 & a_4^3 & a_3^3 & a_2^3 & a_5^3 & a_1^3 & a_4^3  \\
a_2^4 & a_5^4 & a_1^4 & a_4^4 & a_3^4 & a_2^4 & a_5^4 & a_1^4 & a_4^4  \\
a_2^5 & a_5^5 & a_1^5 & a_4^5 & a_3^5 & a_2^5 & a_5^5 & a_1^5 & a_4^5  \\
};
\draw[thick] (mtrx-1-1.north) -| (mtrx-5-1.south west)
                              -- (mtrx-5-1.south);
\draw[thick] (mtrx-1-5.north) -| (mtrx-5-5.south east)
                              -- (mtrx-5-5.south);
\draw [decorate,decoration={brace,amplitude=10pt,raise=2pt},yshift=0pt]
(mtrx-1-1.north west) -- (mtrx-1-5.north east) node [black,midway,yshift=.6cm] {\footnotesize
corresponds to  $(34)\cdot\sigma^2=$\perm{2&5&1&4&3}};
\draw [decorate,decoration={brace,mirror,amplitude=10pt},xshift=3pt]
(mtrx-5-6.south west) -- (mtrx-5-9.south east) node [black,midway,yshift=-.5cm] {\footnotesize
copy first four columns};
\path[draw,strip=orange,shorten <=-1.5mm,shorten >=-1.5mm]
    (mtrx-1-1.center) edge (mtrx-5-5.center)
    (mtrx-1-2.center) edge (mtrx-5-6.center)
    (mtrx-1-3.center) edge (mtrx-5-7.center)
    (mtrx-1-4.center) edge (mtrx-5-8.center)
    (mtrx-1-5.center) edge (mtrx-5-9.center); %
\path[draw,strip=orange,shorten <=-1.5mm,shorten >=-1.5mm]
    (mtrx-5-1.center) edge (mtrx-1-5.center)
    (mtrx-5-2.center) edge (mtrx-1-6.center)
    (mtrx-5-3.center) edge (mtrx-1-7.center)
    (mtrx-5-4.center) edge (mtrx-1-8.center)
    (mtrx-5-5.center) edge (mtrx-1-9.center);
\end{tikzpicture}}\hspace{0cm}
   \xxmini{0.5}{\small \begin{align*}
        N_3 & = \pro{2}{5}{1}{4}{3}+\pro{5}{1}{4}{3}{2}
        + \pro{1}{4}{3}{2}{5}\\
        & +\pro{4}{3}{2}{5}{1} + \pro{3}{2}{5}{1}{4} \\
        & +\pro{3}{4}{1}{5}{2}+\pro{2}{3}{4}{1}{5}
        +\pro{5}{2}{3}{4}{1}\\
        & +\pro{1}{5}{2}{3}{4} +\pro{4}{1}{5}{2}{3}.
    \end{align*}}

    \ \,\hspace{-0.75cm}\xxmini{0.5}{\begin{tikzpicture}[
strip/.style = {
    draw=#1,
    line width=1.25em, opacity=0.2,
    line cap=round ,% only if you like them ...
                            },
                    ]
\matrix (mtrx)  [matrix of math nodes,
                 column sep=0.75em,
                 nodes={text height=2ex,text width=2ex}]
{
a_4^1 & a_1^1 & a_3^1 & a_5^1 & a_2^1 & a_4^1 & a_1^1 & a_3^1 & a_5^1  \\
a_4^2 & a_1^2 & a_3^2 & a_5^2 & a_2^2 & a_4^2 & a_1^2 & a_3^2 & a_5^2  \\
a_4^3 & a_1^3 & a_3^3 & a_5^3 & a_2^3 & a_4^3 & a_1^3 & a_3^3 & a_5^3  \\
a_4^4 & a_1^4 & a_3^4 & a_5^4 & a_2^4 & a_4^4 & a_1^4 & a_3^4 & a_5^4  \\
a_4^5 & a_1^5 & a_3^5 & a_5^5 & a_2^5 & a_4^5 & a_1^5 & a_3^5 & a_5^5  \\
};
\draw[thick] (mtrx-1-1.north) -| (mtrx-5-1.south west)
                              -- (mtrx-5-1.south);
\draw[thick] (mtrx-1-5.north) -| (mtrx-5-5.south east)
                              -- (mtrx-5-5.south);
\draw [decorate,decoration={brace,amplitude=10pt,raise=2pt},yshift=0pt]
(mtrx-1-1.north west) -- (mtrx-1-5.north east) node [black,midway,yshift=.6cm] {\footnotesize
corresponds to  $(34)\cdot\sigma^3=$\perm{4&1&3&5&2}};
\draw [decorate,decoration={brace,mirror,amplitude=10pt},xshift=3pt]
(mtrx-5-6.south west) -- (mtrx-5-9.south east) node [black,midway,yshift=-.5cm] {\footnotesize
copy first four columns};
\path[draw,strip=orange,shorten <=-1.5mm,shorten >=-1.5mm]
    (mtrx-1-1.center) edge (mtrx-5-5.center)
    (mtrx-1-2.center) edge (mtrx-5-6.center)
    (mtrx-1-3.center) edge (mtrx-5-7.center)
    (mtrx-1-4.center) edge (mtrx-5-8.center)
    (mtrx-1-5.center) edge (mtrx-5-9.center); %
\path[draw,strip=orange,shorten <=-1.5mm,shorten >=-1.5mm]
    (mtrx-5-1.center) edge (mtrx-1-5.center)
    (mtrx-5-2.center) edge (mtrx-1-6.center)
    (mtrx-5-3.center) edge (mtrx-1-7.center)
    (mtrx-5-4.center) edge (mtrx-1-8.center)
    (mtrx-5-5.center) edge (mtrx-1-9.center);
\end{tikzpicture}}\hspace{0cm}
   \xxmini{0.5}{\small \begin{align*}
        N_4 & =  \pro{4}{1}{3}{5}{2}+\pro{1}{3}{5}{2}{4}
         +\pro{3}{5}{2}{4}{1}\\
        & +\pro{5}{2}{4}{1}{3}  +\pro{2}{4}{1}{3}{5} \\
        & +\pro{2}{5}{3}{1}{4}+\pro{4}{2}{5}{3}{1} 
        +\pro{1}{4}{2}{5}{3}\\
        & +\pro{3}{1}{4}{2}{5} +\pro{5}{3}{1}{4}{2}.
    \end{align*}}

       \ \,\hspace{-0.75cm}\xxmini{0.5}{\begin{tikzpicture}[
strip/.style = {
    draw=#1,
    line width=1.25em, opacity=0.2,
    line cap=round ,% only if you like them ...
                            },
                    ]
\matrix (mtrx)  [matrix of math nodes,
                 column sep=0.75em,
                 nodes={text height=2ex,text width=2ex}]
{
a_5^1 & a_3^1 & a_2^1 & a_1^1 & a_4^1 & a_5^1 & a_3^1 & a_2^1 & a_1^1  \\
a_5^2 & a_3^2 & a_2^2 & a_1^2 & a_4^2 & a_5^2 & a_3^2 & a_2^2 & a_1^2  \\
a_5^3 & a_3^3 & a_2^3 & a_1^3 & a_4^3 & a_5^3 & a_3^3 & a_2^3 & a_1^3  \\
a_5^4 & a_3^4 & a_2^4 & a_1^4 & a_4^4 & a_5^4 & a_3^4 & a_2^4 & a_1^4  \\
a_5^5 & a_3^5 & a_2^5 & a_1^5 & a_4^5 & a_5^5 & a_3^5 & a_2^5 & a_1^5  \\
};
\draw[thick] (mtrx-1-1.north) -| (mtrx-5-1.south west)
                              -- (mtrx-5-1.south);
\draw[thick] (mtrx-1-5.north) -| (mtrx-5-5.south east)
                              -- (mtrx-5-5.south);
\draw [decorate,decoration={brace,amplitude=10pt,raise=2pt},yshift=0pt]
(mtrx-1-1.north west) -- (mtrx-1-5.north east) node [black,midway,yshift=.6cm] {\footnotesize
corresponds to  $(34)\cdot\sigma^4=$\perm{5&3&2&1&4}};
\draw [decorate,decoration={brace,mirror,amplitude=10pt},xshift=3pt]
(mtrx-5-6.south west) -- (mtrx-5-9.south east) node [black,midway,yshift=-.5cm] {\footnotesize
copy first four columns};
\path[draw,strip=orange,shorten <=-1.5mm,shorten >=-1.5mm]
    (mtrx-1-1.center) edge (mtrx-5-5.center)
    (mtrx-1-2.center) edge (mtrx-5-6.center)
    (mtrx-1-3.center) edge (mtrx-5-7.center)
    (mtrx-1-4.center) edge (mtrx-5-8.center)
    (mtrx-1-5.center) edge (mtrx-5-9.center); %
\path[draw,strip=orange,shorten <=-1.5mm,shorten >=-1.5mm]
    (mtrx-5-1.center) edge (mtrx-1-5.center)
    (mtrx-5-2.center) edge (mtrx-1-6.center)
    (mtrx-5-3.center) edge (mtrx-1-7.center)
    (mtrx-5-4.center) edge (mtrx-1-8.center)
    (mtrx-5-5.center) edge (mtrx-1-9.center);
\end{tikzpicture}}\hspace{0cm}
   \xxmini{0.5}{\small \begin{align*}
        N_5 & = \pro{5}{3}{2}{1}{4}+\pro{3}{2}{1}{4}{5}
         +\pro{2}{1}{4}{5}{3}\\
        & +\pro{1}{4}{5}{3}{2} 
         +\pro{4}{5}{3}{2}{1} \\
        & +\pro{4}{1}{2}{3}{5}+\pro{5}{4}{1}{2}{3} 
        +\pro{3}{5}{4}{1}{2}\\
        & +\pro{2}{3}{5}{4}{1} 
         +\pro{1}{2}{3}{5}{4}.
    \end{align*}}

       \ \,\hspace{-0.75cm}\xxmini{0.5}{\begin{tikzpicture}[
strip/.style = {
    draw=#1,
    line width=1.25em, opacity=0.2,
    line cap=round ,% only if you like them ...
                            },
                    ]
\matrix (mtrx)  [matrix of math nodes,
                 column sep=0.75em,
                 nodes={text height=2ex,text width=2ex}]
{
a_1^1 & a_3^1 & a_2^1 & a_4^1 & a_5^1 & a_1^1 & a_3^1 & a_2^1 & a_4^1  \\
a_1^2 & a_3^2 & a_2^2 & a_4^2 & a_5^2 & a_1^2 & a_3^2 & a_2^2 & a_4^2  \\
a_1^3 & a_3^3 & a_2^3 & a_4^3 & a_5^3 & a_1^3 & a_3^3 & a_2^3 & a_4^3  \\
a_1^4 & a_3^4 & a_2^4 & a_4^4 & a_5^4 & a_1^4 & a_3^4 & a_2^4 & a_4^4  \\
a_1^5 & a_3^5 & a_2^5 & a_4^5 & a_5^5 & a_1^5 & a_3^5 & a_2^5 & a_4^5  \\
};
\draw[thick] (mtrx-1-1.north) -| (mtrx-5-1.south west)
                              -- (mtrx-5-1.south);
\draw[thick] (mtrx-1-5.north) -| (mtrx-5-5.south east)
                              -- (mtrx-5-5.south);
\draw [decorate,decoration={brace,amplitude=10pt,raise=2pt},yshift=0pt]
(mtrx-1-1.north west) -- (mtrx-1-5.north east) node [black,midway,yshift=.6cm] {\footnotesize
corresponds to  $(34)\cdot\tau=$\perm{1&3&2&4&5}};
\draw [decorate,decoration={brace,mirror,amplitude=10pt},xshift=3pt]
(mtrx-5-6.south west) -- (mtrx-5-9.south east) node [black,midway,yshift=-.5cm] {\footnotesize
copy first four columns};
\path[draw,strip=orange,shorten <=-1.5mm,shorten >=-1.5mm]
    (mtrx-1-1.center) edge (mtrx-5-5.center)
    (mtrx-1-2.center) edge (mtrx-5-6.center)
    (mtrx-1-3.center) edge (mtrx-5-7.center)
    (mtrx-1-4.center) edge (mtrx-5-8.center)
    (mtrx-1-5.center) edge (mtrx-5-9.center); %
\path[draw,strip=orange,shorten <=-1.5mm,shorten >=-1.5mm]
    (mtrx-5-1.center) edge (mtrx-1-5.center)
    (mtrx-5-2.center) edge (mtrx-1-6.center)
    (mtrx-5-3.center) edge (mtrx-1-7.center)
    (mtrx-5-4.center) edge (mtrx-1-8.center)
    (mtrx-5-5.center) edge (mtrx-1-9.center);
\end{tikzpicture}}\hspace{0cm}
   \xxmini{0.5}{\small \begin{align*}
        N_6 & = \pro{1}{3}{2}{4}{5}+\pro{3}{2}{4}{5}{1} 
        +\pro{2}{4}{5}{1}{3}\\
        & +\pro{4}{5}{1}{3}{2}
        +\pro{5}{1}{3}{2}{4}\\
        & +\pro{5}{4}{2}{3}{1}+\pro{1}{5}{4}{2}{3}
         +\pro{3}{1}{5}{4}{2}\\
        & +\pro{2}{3}{1}{5}{4}
         +\pro{4}{2}{3}{1}{5}.
    \end{align*}}

Similarly to the previous set of matrices, when putting together all these augmented matrices--which correspond to all the odd permutations--they connect in such a way that the last column of one is the first column of the next one. Here are the first two matrices together:

\medskip

 \ \,\hspace{-0.5cm}\xxmini{0.5}{\begin{tikzpicture}[
strip/.style = {
    draw=#1,
    line width=1.25em, opacity=0.2,
    line cap=round ,% only if you like them ...
                            },
                    ]
\matrix (mtrx)  [matrix of math nodes,
                 column sep=1em,
                 nodes={text height=2ex,text width=2ex}]
{
a_1^1 & a_2^1 & a_4^1 & a_3^1 & a_5^1 & a_1^1 & a_2^1 & a_4^1 &
a_3^1 & a_4^1 & a_5^1 & a_2^1 & a_1^1 & a_3^1 & a_4^1 & a_5^1 & a_2^1 \\
a_1^2 & a_2^2 & a_4^2 & a_3^2 & a_5^2 & a_1^2 & a_2^2 & a_4^2 &
a_3^2 & a_4^2 & a_5^2 & a_2^2 & a_1^2 & a_3^2 & a_4^2 & a_5^2 & a_2^2   \\
a_1^3 & a_2^3 & a_4^3 & a_3^3 & a_5^3 & a_1^3 & a_2^3 & a_4^3 &
a_3^3 & a_4^3 & a_5^3 & a_2^3 & a_1^3 & a_3^3 & a_4^3 & a_5^3 & a_2^3   \\
a_1^4 & a_2^4 & a_4^4 & a_3^4 & a_5^4 & a_1^4 & a_2^4 & a_4^4 &
a_3^4 & a_4^4 & a_5^4 & a_2^4 & a_1^4 & a_3^4 & a_4^4 & a_5^4 & a_2^4   \\
a_1^5 & a_2^5 & a_4^5 & a_3^5 & a_5^5 & a_1^5 & a_2^5 & a_4^5 &
a_3^5 & a_4^5 & a_5^5 & a_2^5 & a_1^5 & a_3^5 & a_4^5 & a_5^5 & a_2^5   \\
};
\draw[thick] (mtrx-1-1.north) -| (mtrx-5-1.south west)
                              -- (mtrx-5-1.south);
\draw[thick] (mtrx-1-5.north) -| (mtrx-5-5.south east)
                              -- (mtrx-5-5.south);
\draw [decorate,decoration={brace,amplitude=10pt,raise=2pt},yshift=0pt]
(mtrx-1-1.north west) -- (mtrx-1-5.north east) node [black,midway,yshift=.6cm] {\footnotesize
corresponds to  $(34)\cdot e=$\perm{1&2&4&3&5}}; 
\path[draw,strip=orange,shorten <=-1.5mm,shorten >=-1.5mm]
    (mtrx-1-1.center) edge (mtrx-5-5.center)
    (mtrx-1-2.center) edge (mtrx-5-6.center)
    (mtrx-1-3.center) edge (mtrx-5-7.center)
    (mtrx-1-4.center) edge (mtrx-5-8.center)
    (mtrx-1-5.center) edge (mtrx-5-9.center); %
\path[draw,strip=orange,shorten <=-1.5mm,shorten >=-1.5mm]
    (mtrx-5-1.center) edge (mtrx-1-5.center)
    (mtrx-5-2.center) edge (mtrx-1-6.center)
    (mtrx-5-3.center) edge (mtrx-1-7.center)
    (mtrx-5-4.center) edge (mtrx-1-8.center)
    (mtrx-5-5.center) edge (mtrx-1-9.center);
\draw [decorate,decoration={brace,amplitude=10pt,raise=2pt},yshift=0pt]
(mtrx-1-9.north west) -- (mtrx-1-13.north east) node [black,midway,yshift=.6cm] {\footnotesize
corresponds to  $(34)\cdot\sigma=$\perm{3&4&5&2&1}}; 
\path[draw,strip=orange,shorten <=-1.5mm,shorten >=-1.5mm]
    (mtrx-1-9.center) edge (mtrx-5-13.center)
    (mtrx-1-10.center) edge (mtrx-5-14.center)
    (mtrx-1-11.center) edge (mtrx-5-15.center)
    (mtrx-1-12.center) edge (mtrx-5-16.center)
    (mtrx-1-13.center) edge (mtrx-5-17.center); %
\path[draw,strip=orange,shorten <=-1.5mm,shorten >=-1.5mm]
    (mtrx-5-9.center) edge (mtrx-1-13.center)
    (mtrx-5-10.center) edge (mtrx-1-14.center)
    (mtrx-5-11.center) edge (mtrx-1-15.center)
    (mtrx-5-12.center) edge (mtrx-1-16.center)
    (mtrx-5-13.center) edge (mtrx-1-17.center);
\end{tikzpicture}}

Here is all of this set of matrices sewn together into our second ``quilt":

\bigskip

\ \,\hspace{-0.65cm}\includegraphics[width=0.75\paperwidth]{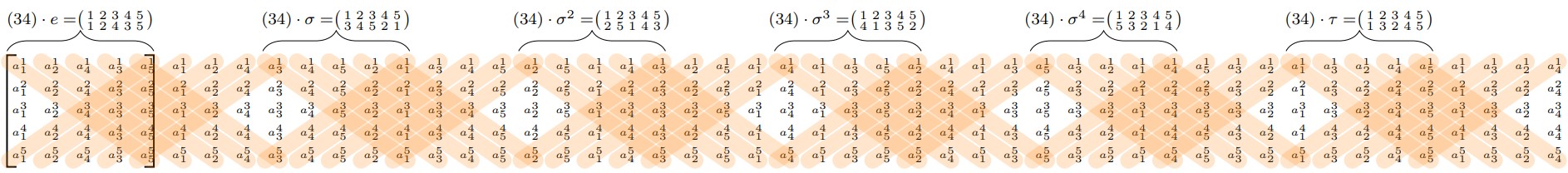}

\section{Larger Determinants}\label{LargerDeterm}
We have observed there to be a repeating pattern as we increased the size of a $n \times n$ matrix. As we set up a $6 \times 6$ matrix, we noticed that the products of the diagonals were organized in the same manner as a $2 \times 2$ matrix. We continued again to a $7 \times 7$ matrix and noticed a similar structure to a $3 \times 3$ matrix. Repetition of this has led us to find that every $n \times n$ matrix falls into one of four patterns, where $n=4k+2$, $n=4k+3$, $n=4k+4$, or $n=4k+5$, such that $k \in \mathbb{Z}_{\ge 0}$.

 \ \,\hspace{-0.5cm} \begin{minipage}[t]{\textwidth}

\centering
\begin{tikzpicture}[framed,
strip/.style = {
    draw=#1,%color
    line width=1em, opacity=0.2,
    line cap=round ,% only if you like them ...
                            },
                    ]
\matrix (mtrx)  [matrix of math nodes,
                 column sep=1em,
                 nodes={text height=2ex,text width=1ex}]
{
\, & \, &\, &\, &\, &\, &\, \\
\, & \, &\, &\, & \,&\cdots &\, \\
\, & \, &\, &\, &\, &\, &\, \\[-2mm]
};
\path[draw,strip=cyan,shorten <=-1mm]
    (mtrx-1-1.center) edge (mtrx-3-3.center)
    (mtrx-1-3.center) edge (mtrx-3-5.center)
    (mtrx-3-2.center) edge (mtrx-1-4.center)
    (mtrx-3-4.center) --   (mtrx-1-6.center); 
\path[draw,strip=orange,shorten <=-1mm]
    (mtrx-1-3.center) edge (mtrx-3-1.center)
    (mtrx-1-5.center) edge (mtrx-3-3.center)
    (mtrx-1-4.center) edge (mtrx-3-6.center)
    (mtrx-1-2.center) --   (mtrx-3-4.center);
\draw [decorate,decoration={brace,amplitude=10pt,raise=2pt},yshift=0pt]
(mtrx-1-1.north west) -- (mtrx-1-6.north east) node [black,midway,yshift=.6cm] {\footnotesize Pattern when $n=4k+2$}; 

\end{tikzpicture}   
\begin{tikzpicture}[framed,
strip/.style = {
    draw=#1,%color
    line width=1em, opacity=0.2,
    line cap=round ,% only if you like them ...
                            },
                    ]
\matrix (mtrx)  [matrix of math nodes,
                 column sep=1em,
                 nodes={text height=2ex,text width=1ex}]
{
\, & \, &\, &\, &\, &\, &\,\\
\, & \, &\, &\, &\, &\cdots &\, \\
\, & \, &\, &\, &\, &\, &\, \\[-2mm]
};
\path[draw,strip=cyan,shorten <=-1mm]
    (mtrx-1-1.center) edge (mtrx-3-3.center)
    (mtrx-1-2.center) edge (mtrx-3-4.center)
    (mtrx-1-3.center) edge (mtrx-3-5.center)
    (mtrx-1-4.center) -- (mtrx-3-6.center); 
\path[draw,strip=orange,shorten <=-1mm]
    (mtrx-1-3.center) edge (mtrx-3-1.center)
    (mtrx-1-4.center) edge (mtrx-3-2.center)
    (mtrx-1-5.center) edge (mtrx-3-3.center)
    (mtrx-1-6.center) --   (mtrx-3-4.center); 
\draw [decorate,decoration={brace,amplitude=10pt,raise=2pt},yshift=0pt]
(mtrx-1-1.north west) -- (mtrx-1-6.north east) node [black,midway,yshift=.6cm] {\footnotesize Pattern when $n=4k+3$}; 
\end{tikzpicture}
\end{minipage} 

\ \, \vspace{-3mm}\ \,

\begin{minipage}[t]{\textwidth}
 \centering
      \begin{tikzpicture}[framed,
strip/.style = {
    draw=#1,%color
    line width=1em, opacity=0.2,
    line cap=round ,% only if you like them ...
                            },
                    ]
\matrix (mtrx)  [matrix of math nodes,
                 column sep=1em,
                 nodes={text height=2ex,text width=1ex}]
{
\, & \, &\, &\, &\, &\, &\, \\
\, & \, &\, &\, & \,&\cdots &\, \\
\, & \, &\, &\, &\, &\, &\, \\[-2mm]
};
\path[draw,strip=cyan,shorten <=-1mm]
    (mtrx-1-1.center) edge (mtrx-3-3.center)
    (mtrx-3-1.center) edge (mtrx-1-3.center)
    (mtrx-1-3.center) edge (mtrx-3-5.center)
    (mtrx-3-3.center) --   (mtrx-1-5.center); 
\path[draw,strip=orange,shorten <=-1mm]
    (mtrx-1-2.center) edge (mtrx-3-4.center)
    (mtrx-1-4.center) edge (mtrx-3-2.center)
    (mtrx-1-4.center) edge (mtrx-3-6.center)
    (mtrx-1-6.center) --   (mtrx-3-4.center);
\draw [decorate,decoration={brace,amplitude=10pt,raise=2pt},yshift=0pt]
(mtrx-1-1.north west) -- (mtrx-1-6.north east) node [black,midway,yshift=.6cm] {\footnotesize Pattern when $n=4k+4$}; 
\end{tikzpicture}
      \begin{tikzpicture}[framed,
strip/.style = {
    draw=#1,%color
    line width=1em, opacity=0.2,
    line cap=round ,% only if you like them ...
                            }
                    ]
\matrix (mtrx)  [matrix of math nodes,
                 column sep=1em,
                 nodes={text height=2ex,text width=1ex}]
{
\, & \, &\, &\, \, & \, &\, &\, &\,\\
\, & \, &\, &\cdots\, & \, &\, &\, &\cdots \\
\, & \, &\, &\, \, & \, &\, &\, &\,\\[-2mm]
};
\path[draw,strip=cyan,shorten <=-1mm]
    (mtrx-1-1.center) edge (mtrx-3-3.center)
    (mtrx-3-1.center) edge (mtrx-1-3.center)
    (mtrx-1-2.center) edge (mtrx-3-4.center)
    (mtrx-1-4.center) edge (mtrx-3-2.center); 
\path[draw,strip=orange,shorten <=-1mm]
    (mtrx-1-5.center) edge (mtrx-3-7.center)
    (mtrx-3-5.center) edge (mtrx-1-7.center)
    (mtrx-1-6.center) edge (mtrx-3-8.center)
    (mtrx-1-8.center) edge (mtrx-3-6.center);
\draw [decorate,decoration={brace,amplitude=10pt,raise=2pt},yshift=0pt]
(mtrx-1-1.north west) -- (mtrx-1-8.north east) node [black,midway,yshift=.6cm] {\footnotesize Pattern when $n=4k+5$}; 
\end{tikzpicture} 
\end{minipage}

\section{Comparison with Existing Methods}\label{ComparisonExisting}

Our proposed method constructs a $4 \times 19$ matrix by repeating and reordering the columns of the original matrix. Diagonals are taken directly from this matrix, with signs assigned based on the starting column of each diagonal. This structure captures all $24$ terms required for the determinant, similar in spirit to the Sarrus rule for $3 \times 3$ matrices. The approach is straightforward as it does not
%, visual, and easy to follow, 
require minor expansions, row or column operations, and multiple matrices.

The method introduced by Salihu and Gjonbalaj (2017, see~\cite{Armend4x4}) differs by using three separate row-permuted matrices to generate the necessary diagonals. Each matrix contributes eight diagonals to cover all terms of the determinant. While their approach is systematic and avoids cofactor expansion, it involves multiple matrix constructions and shifts the user's attention across different layouts. In contrast, our method achieves the same goal using a single fixed structure.
% making the process more streamlined and visually cohesive, especially for manual or instructional use.

Sobamowo’s (2016, see~\cite{Sobamowo}) approach aims to generalize Sarrus' rule to $n \times n$ matrices by cyclically permuting columns and applying a structured diagonal scheme. This method provides scalability, but increases procedural complexity and has a more complex visual structure. Our method improves usability by eliminating column permutations and embedding the entire process into a fixed-width matrix, making the sign pattern and diagonal traversal more comprehensible and replicable without reference to external rules.

Salinas-Hernández et al. (2021, see~\cite{Salinas}) propose a method that extends the original matrix to $6 \times 4$ by duplicating specific rows and applies a “ladder” diagonal multiplication to each term in the first row. This method retains ties to co-factor expansion and effectively reduces computational effort. However, it still requires separate evaluations per top-row element and more structural manipulation of the matrix. Our method avoids both row duplication and per-element diagonal sets, enabling all required diagonals to be extracted from a single pass through a column-augmented matrix.
%, reducing both mental load and mechanical steps.

The method presented by Oliveira (2013, see~\cite{Oliveira}) extends Sarrus’ rule to $5 \times 5$ matrices by constructing a $5 \times 13$ matrix formed by repeating and shifting the original matrix’s columns. Diagonal products are then extracted using a consistent pattern, with alternating signs assigned to produce the full determinant. While the method preserves the visual and intuitive strengths of Sarrus’ approach, it is specifically developed for the $5 \times 5$ case and requires more column manipulation. In contrast, our method in the $5\times 5$ case, connects all the matrices together in a ``quilt-pattern" via specific permutations.

Arshon (1935, see~\cite{Arshon}) introduces a generalization of the Sarrus rule for $n \times n$ matrices by decomposing the determinant into a sum of $(n-1)!/2$ Sarrus-style diagrams, each derived from a distinct permutation of the columns. For each permutation, a $3 \times 3$ window is evaluated using traditional Sarrus diagonals, and the results are aggregated with appropriate sign assignments to reconstruct the full $n!$-term determinant expansion. While this method offers a compelling graphical and combinatorial interpretation of the determinant, it requires managing multiple column permutations and tracking sign contributions for each diagram. In contrast, our proposed method consolidates all necessary diagonals into a single $4 \times 19$ matrix.

In comparison with all these works, our method provides patterns for larger determinants and how they behave when $n=4k+i, 
 i=0,1,2,3$. 

%Add in the two papers we discussed in class, cite and compare them in there
% https://rsdjournal.org/index.php/rsd/article/view/40121/32822 - Portuguese
% https://www.mathnet.ru/links/e6f0077a8cd7a979588f64dd19bb9976/sm6404.pdf - Russian

%\subsection*{Summary of Differences}

%\begin{itemize}
%    \item \textbf{Single-structure layout:} Unlike others, our method uses only one matrix for computation.
%    \item \textbf{Intuitive diagonal logic:} Sign rules are tied to clear column indices, not permutations or recursive logic.
%    \item \textbf{Manual and visual ease:} It is well suited for teaching, calculation by hand, and algorithmic implementation.
%\end{itemize}

 \bibliography{Sarrus}

\end{document}